\documentclass{article}
\usepackage[utf8]{inputenc}
\usepackage[english]{babel}
\usepackage{amsmath,amsfonts,amssymb,amsthm}
\usepackage{indentfirst}

\newtheorem{definition}{Definition}

\newtheorem{proposition}{Proposition}
\newtheorem{remark}{Remark}

\numberwithin{equation}{section}
\newcommand{\R}{\ensuremath{\mathbb{R}}}

\newcommand{\Z}{\ensuremath{\mathbb{Z}}}

\newcommand{\ep}{\varepsilon}

\newcommand{\de}{\delta}
\newcommand{\De}{\Delta}
\newcommand{\dist}{\textup{dist}}

\newcommand{\Image}{\textup{Im}}
\newcommand{\conv}{\textup{conv}}
\title{Multidimensional $C^0$ transversality}
\author{Aleksey A. Petrov, Sergei Yu. Pilyugin}

\begin{document}

\maketitle

\begin{abstract}

In 1994, Sakai introduced the property of $C^0$ transversality for two smooth curves 
in a two-dimensional manifold. This property was related to various shadowing properties
of dynamical systems. In this short note, we generalize this property to arbitrary
continuous mappings of topological spaces into topological manifolds. We prove 
a sufficient condition for the $C^0$ transversality of two submanifolds of a topological manifold
and a necessary condition of $C^0$ transversality for mappings of metric spaces into $\R^n$.
Bibliography: 7 titles.
\medskip

Keywords: continuous mapping, manifold, topological space, dimension
\medskip

AMS Subject Classification: 54C05

\end{abstract}

\section{Introduction}

In the paper \cite{Sak}, Sakai formulated the so-called condition of $C^0$ transversality for two
smooth curves in a two-dimensional manifold and gave a sketch of the proof of the following
theorem: For an Axiom A diffeomorphism $f$ of a two-dimensional manifold, the following two statements
are equivalent:
\begin{itemize}
\item $f$ has the shadowing property;
\item one-dimensional stable and unstable manifolds of hyperbolic nonwandering points
are  $C^0$-transverse.
\end{itemize}

(One can find necessary definitions related to shadowing theory for dynamical systems
in the book \cite{PilSpDS}.)

Later, a complete proof of the above-mentioned theorem was given in \cite{PilSak}
(it was also shown in \cite{PilSak} that a two-dimensional Axiom A diffeomorphism
has the inverse shadowing property with respect to continuous methods if and only
if it satisfies the above $C^0$ transversality condition).

In this paper, we generalize the notion of $C^0$ transversality to arbitrary continuous
mappings of topological spaces into $\R^n$ (or, more general, into topological
manifolds). In our opinion, this notion may be of interest not only for specialists
in dynamical systems.

The structure of the paper is as follows. In Sec. 2, we define multidimensional $C^0$ transversality.
Section 3 contains a sufficient condition for the $C^0$ transversality of two submanifolds of a topological manifold.
We also show that in the case of two smooth curves in a two-dimensional manifold, our
definition coincides with that given in \cite{Sak}. Section 4 is devoted to a necessary
condition of $C^0$ transversality for mappings of metric spaces into $\R^n$.

Everywhere in this paper, a manifold is a Hausdorff paracompact topological space
with a fixed structure of a topological manifold.

\section{Basic definitions}

In this section, we define $C^0$ transversality for two continuous mappings of
topological spaces into a manifold.

Let $(M, \dist)$ be a topological manifold with a fixed metric $\dist$ and let $A$ be a topological space.

On the space of continuous mappings of $A$ into $M$ (denoted $C(A, M)$) we consider the $C^0$ uniform metric, 
$$
|f_1,f_2|_{C^0}=\sup_{x\in A }(\dist(f_1(x),f_2(x)))
$$
for $f_1,f_2\in C(A, M)$.

\begin{definition}
Let $A$ and $B$ be topological spaces and let $h_1\colon A\to M$ and $h_2\colon B\to M$
be continuous mappings. If $U_A\subset A$ and $U_B\subset B$ are arbitrary subsets and
$\delta>0$, we say that the intersection $h_1(U_A)\cap h_2(U_B)$ is $\delta$-essential if
$$
\widetilde{h_1}(U_A) \cap \widetilde{h_2}(U_B)\neq \varnothing
$$
for any continuous mappings 
$$
\widetilde{h_1}\colon A\to M
\quad\mbox{and}\quad
\widetilde{h_2}\colon B\to M
$$
such that $|h_1,\widetilde{h_1}|_{C^0}\leq \de$ and $|h_2,\widetilde{h_2}|_{C^0}\leq \de$.
\end{definition}

\begin{remark} In the classical book \cite{GW}, the authors condidered a similar property
in the particular case where $B$ is a point and $\widetilde{h_2}=h_2$ 
(the corresponding point $h_2(B)$ was called a stable
value of the mapping $h_1$). They proved the following results concerning continuous 
mappings into the $n$-dimensional cube: If $\dim(A)<n$, then 
any value of $h_1$ is unstable,
and if $\dim(A)\geq n$, then there exists a mapping $h_1$ having at least one
stable value (compare with Proposition 4 below).
\end{remark}

\begin{definition} 
Let again $A$ and $B$ be topological spaces and let $h_1\colon A\to M$ and $h_2\colon B\to M$
be continuous mappings. Assume that $h_1(a)=h_2(b)$ for some points $a\in A$ and $b\in B$.
We say that the mappings $h_1$ and $h_2$ are $C^0$ transverse at the pair $(a,b)$ if for any
open sets $U(a)\subset A$ and $U(b)\subset B$ such that $a\in U(a)$ and $b\in U(b)$
there exists a $\delta>0$ such that the intersection $h_1(U(a))\cap h_2(U(b))$
is $\delta$-essential.
\end{definition} 

Finally, we define the $C^0$ transversality of two mappings.
\begin{definition}
Let $A$ and $B$ be topological spaces and let $h_1\colon A\to M$ and $h_2\colon B\to M$
be continuous mappings. We say that $h_1$ and $h_2$ are $C^0$ transverse if
$h_1$ and $h_2$ are $C^0$ transverse at any pair
$(a,b)\in A\times B$ such that $h_1(a)=h_2(b)$.
\end{definition}

Clearly, the given definition of $C^0$ transversality for two mappings does not depend
on the metric fixed on $M$. In the particular case where $A, B\subset M$, it is convenient 
for us to give one more definition.

\begin{definition} 
Let $A\subset M, B\subset M$. We say that $A$ and $B$ are $C^0$ transverse at a point
$p\in A\cap B$  $A$ if the embeddings $i_A\colon A \to M$, $i_A(x)=x$ and $i_B\colon B \to M$, $i_B(x)=x$, 
are $C^0$ transverse at the pair $(p,p)$.
\end{definition} 

Note that if $M$ is a smooth manifold and $A$ and $B$ are smooth submanifolds $M$ without boundary such that 
$A$ and $B$ are transverse at a point $p\in A\cap B$ in the sense of the standard definition from differential
topology (i.e., $T_pA+T_pB=T_pM$), then $A$ and $B$ are $C^0$ transverse at $p$. In particular, this
follows from our Proposition~2 proved in Sec. 3.

The following statement is a corollary of the possibility of approximation of a continuous mapping by a smooth
one and from the standard smooth transversality theorem (see, for example, \cite{Hirsch}). 

\begin{proposition}
Let  $A$, $B$, and $M$ be connected smooth manifolds and let $f\colon A \to M$ and $g\colon B \to M$ 
be continuous mappings for which there exists a pair of points $(a,b)\in A\times B$
such that $f$ and $g$ are $C^0$ transverse at $(a,b)$. Then  $\dim(A)+\dim(B)\geq \dim(M)$.
\end{proposition}

Note that an analog of this statement holds for nonsmooth manifolds as well
(see Proposition~4 in Sec.~4).

\section{Sufficient conditions}

In applications to dynamical systems, we are mostly interested in intersections of
two submanifolds of a manifold. In this section, we prove a simple sufficient condition
for the  $C^0$ transversality of two topological submanifolds $A$ and $B$ of a manifold $M$
(in fact, we give a condition under which the intersection of two open subsets of the
submanifolds is $\delta$-essential for small $\delta$).
Since the definition of $C^0$ transversality is local, without loss of
generality we may assume that $M=\R^n$.

In the space $\R^n$, we consider the norm 
$$
|(a_1,\dots, a_n)|=\max\{ |a_1|, \dots, |a_n|\},
$$
where $a_i$ is the $i$th coordinate of a vector.

Let $A,B\subset \R^n$ be topological submanifolds and let $p\in A\cap B$.  
Let  $k=\dim (A)$ and $l=\dim (B)$.

Denote by $i_A\colon A\to \R^n$ and $i_B\colon B \to \R^n$ the embeddings of $A$ and $B$, respectively.

Since $A$ is a submanifold of $\R^n$, there exists a neighborhood $U\subset \R^n$ of $p$ 
and a homeomorphism
\begin{equation}
\phi\colon U \to \R^n
\end{equation}
such that $\phi(p)=0$ and 
$$
\phi(A\cap U)=\R^k\times \{0\}^{n-k}=\{ (a_1,\dots, a_k, 0, \dots, 0)\in \R^n \mid a_1,\dots, a_k\in \R  \}.
$$

Let $V_A\subset A$ be a neighborhood of $p$ in $A$ such that $Cl(V_A) \subset U$ and
$V_B\subset B$ be a neighborhood of $p$ in $B$ such that $Cl(V_B) \subset U$.

Fix a $\nu>0$ and set
\begin{equation}
J=D^n=[-1,1]^{n}
\end{equation}
and
\begin{equation}
J_A=[-1-\nu,1+\nu ]^{k}\times \{0\}^{n-k}.
\end{equation}

Without loss of generality, we may assume that $J_A\subset \phi(V_A)$.

%Не умаляя общности, можем считать, что $J_A\subset \phi(V_A)$ (в противном случае, мы можем рассмотреть композицию   $\phi$
%с растягивающей гомотетией с центром в $0$).

Since $J\setminus J_A$ is homotopically equivalent to $S^{n-k-1}=\partial D^{n-k}$,
\begin{equation}
\widetilde{H}_{n-k-1}(J\setminus J_A)=\Z,
\end{equation}
where $\widetilde{H}_{n-k-1}$ is the (n-k-1)st reduced homology group.

Now we formulate \textbf{Condition $T(V_A,V_B):$}\  
There exists a continuous mapping 
$$
f\colon D^{n-k}\to V_B
$$
such that
\begin{equation}
\phi(f(D^{n-k}))\subset J,
\end{equation}
\begin{equation}
\phi(f(S^{n-k-1}))\subset J\setminus J_A,
\end{equation}
and the induced homomorphism of the reduced homology groups,
\begin{equation}
(\phi\circ f)_{*}\colon \widetilde{H}_{n-k-1}(S^{n-k-1})\to \widetilde{H}_{n-k-1}(J\setminus J_A),
\end{equation}
is nontrivial.

\begin{proposition}
If Condition $T(V_A,V_B)$ holds, then there exists a number $\delta>0$ such that the intersection $i_A(V_A)\cap i_B(V_B)$ 
is $\delta$-essential.
\end{proposition}

\begin{remark} This statement is, in a sense, folklore, but since we cannot find a proper
reference, we include its proof.
\end{remark}

\textit{Proof. } To simplify presentation, we identify $V_A$ with $\phi(V_A)$ and $V_B$ with $\phi(V_B)$.
It is enough to show that there exists a $\delta>0$ such that for any continuous mappings
\begin{equation}
h_A\colon V_A \to \R^n,
\end{equation}
and 
\begin{equation}
h_B \colon V_B \to \R^n,
\end{equation}
with
\begin{equation}
\label{smallA}
|h_A(x)-x|<\delta,\quad x\in V_A,
\end{equation}
\begin{equation}
\label{smallB}
|h_B(x)-x|<\delta,\quad x\in V_B,
\end{equation}
\begin{equation}
h_A(V_A)\cap h_B(V_B)\neq \varnothing.
\end{equation}

During the proof, we decrease $\delta>0$ several times.

Let $\delta>0$ be so small that if $h_A$ and $h_B$ satisfy (\ref{smallA}) and (\ref{smallB}), then
\begin{equation}
\label{13}
h_A(J_A)\cap h_B(f(S^{n-k-1}))=\varnothing,
\end{equation}
\begin{equation}
\label{14}
h_A(\partial J_A)\cap J = \varnothing,
\end{equation}
and
\begin{equation}
\label{15}
h_B(f(D^{n-k})) \cap \{ (a_1,\dots, a_k, 0, \dots, 0)\mid \ |(a_1,\dots, a_k, 0, \dots, 0)|\geq 1+\frac{\nu}{2} \}=\varnothing.
\end{equation}

Denote by $e_1$ a generator of the group $\widetilde{H}_{n-k-1}(S^{n-k-1})$ and by $e_2$ 
a generator of the group $\widetilde{H}_{n-k-1}(J\setminus J_A)$, respectively. Then we can write 
$$
f_{*}(e_1)=\kappa e_2,
$$ 
where
\begin{equation}
f_{*}\colon \widetilde{H}_{n-k-1}(S^{n-k-1}) \to \widetilde{H}_{n-k-1}(J\setminus J_A)
\end{equation}
is the induced homomorphism of the reduced homology groups and $\kappa\in \Z$. 

By Condition $T(V_A,V_B)$ ,
\begin{equation}
\kappa\neq 0.
\end{equation}

Let $h_A$ and $h_B$ satisfy conditions (\ref{smallA}) and (\ref{smallB}). 
To complete the proof, we show that
\begin{equation}
h_A(J_A)\cap h_B(f(D^{n-k}))\neq \varnothing.
\end{equation} 
For this purpose, consider the mapping 
$$
g\colon \partial(I^{k+1})\to \R^n,
$$
where  $I=[-1-\nu,1+\nu]$, defined as follows:
$$
g(a_1,\dots, a_k,1+\nu)= (a_1,\dots, a_k, 0, \dots, 0)
$$ 
if $ |a_i|< 1+\nu, \ i=1,\dots, k,$ and
$$
g(a_1,\dots, a_k,t)=(a_1,\dots,a_k, 3(t-1-\nu),0, \dots, 0)
$$
if one of the values $a_i$ satisfies the relation $|a_i|=1+\nu$.

Since $\Image f \cap \Image g=\varnothing$, the linking coefficient $K(f,g)$ is defined. 
Clearly, the absolute value of the linking coefficient of $f$ and $g$ coincides with
the absolute value of $\kappa$: 
$$
|K(f,g)|=|\kappa|.
$$

Reducing $\delta>0$, if necessary, we may assume that for any two continuous mappings
$$
\widetilde{f}\colon S^{n-k-1}\to \R^n\quad\mbox{and}\quad
\widetilde{g}\colon \partial(I^{k+1})\to \R^n
$$
such that
\begin{equation}
\label{3star}
|\widetilde{f},f|_{C^0}<\delta \quad\mbox{and}\quad
|\widetilde{g},g|_{C^0}<\delta,
\end{equation}
the equality $K(\widetilde{f},\widetilde{g})=K(f,g)$ holds. 

Set $\widetilde{f}=h_B \circ f$ and define a mapping
$\widetilde{g}\colon \partial(I^{k+1}) \to \R^n$ as follows:
$$
\widetilde{g}(a_1,\dots, a_k, 1+\nu)=h_A(a_1,\dots, a_k, 0, \dots, 0)
$$
for all $a_1,\dots, a_k$ such that $|a_i|< 1+\nu$, and
$$
\widetilde{g}(a_1,\dots, a_k, t)=h_A(a_1,\dots, a_k, 0, \dots, 0)+3(t-1-\nu)v_{k+1}
$$
if one of the values $|a_i|$ equals $1+\nu$. Here $v_{k+1}$ is the unit $(k+1)$st basic vector in $\R^n$ 
(i.e., its $(k+1)$st coordinate equals $1$ while the remaining coordinates equal 0).

By construction, $\widetilde{f}$ and $\widetilde{g}$ satisfy conditions (\ref{3star}). Thus, 
$$
|K(\widetilde{f},\widetilde{g})|=|K(f,g)|=|\kappa|\neq 0.
$$
It follows that 
$$
\Image \widetilde{g}\cap \widetilde{f}(D^{n-k}) \neq \varnothing;
$$
otherwise, $K(\widetilde{f},\widetilde{g})=0$.

Finally, conditions (\ref{14}) and (\ref{15}) imply that 
$$
\Image \widetilde{g}\cap \widetilde{f}(D^{n-k})\subset  h_A(J_A)\cap\widetilde{f}(D^{n-k})= h_A(J_A)\cap h_B(f(D^{n-k})),
$$
%и поскольку выполнено включение $$\Image h_A \cap \Image h_b \supset \widetilde{f}(D^{n-k}) \cap \Image h_A,$$
%то $$\Image h_A \cap \Image h_b\neq \varnothing,$$
and we conclude that 
$$h_A(J_A)\cap h_B(f(D^{n-k}))\neq \varnothing,
$$ 
as required.
$\Box $

\begin{remark} It is not clear whether it is possible to find reasonable necessary conditions of the same kind.
Consider, for example, $M=\R^2$ with coordinates $(x,y)$ and let $A$ be the axis $\{(x,0)\}$ while $B$
is the graph of the function
$$
y(x)=
\Big{\{}
\begin{array}{cc}
0,&\quad x\leq 0,\\
\exp(-1/x^2)\sin(1/x),&\quad x>0.\\
\end{array}
$$
Clearly, $A$ and $B$ are $C^0$ transverse at the origin.
\end{remark}

Now we recall the definition of $C^0$ transversality formulated in  \cite{Sak}. 
Let $K$ and $L$ be two $C^1$ curves in a two-dimesional manifold $M$ (i.e., one-dimensional boundaryless $C^1$-submanifolds of $M$). 
Consider a point $p\in K\cap L$ and let $D_{a}$ be the open disk of small radius $a$ in $M$ centered at $p$. 
If $a$ is small enough, then $K'$, the component of $K\cap D_a$ containing $p$, separates the disk $D_{a}$ into two open
components $D_a^+$ and $D_a^-$ (each of them is homeomorphic to an open disk). We say that the curves $K$ and $L$ are 
$C^0$ transverse at $p$ if for any $a>0$, the component $L'$ of the intersection $L\cap D_a$ containing the point $p$
intersects both $D_a^+$ and $D_a^-$. 

Let us show that this definition is equivalent to our Definition 4.
Proposition 2 implies that our Definition 4 is a corollary of the definition given in \cite{Sak}. 
Let us prove the converse. Assume that two $C^1$ curves $K$ and $L$ are $C^0$ transverse at $p$ in the sense of
our Definition 4 but do not satisfy the definition given in \cite{Sak}. Then, if $a$ is small enough,
$L'$ belongs either to
$D_a^+\cup \left( K' \cap D_a \right)$ either to $D_a^-\cup \left( K' \cap D_a \right)$
(for definiteness, assume that $L'\subset D_a^+\cup \left( K' \cap D_a \right)$).
In this case, applying a continuous mapping that is arbitraily close to the identical mapping of $\R^2$,
we can ``squeeze out'' $L'$ into $D_a^+$. In this case, the image of $L'$ would not intersect
$K'$, which contradicts the existence of $\delta$ in Definition 2.

\begin{remark}
It was noted in \cite{PilSakTar}, Remark 2.2 that it is possible to reformulate the definition
of $C^0$ transversality of two smooth curves given in \cite{Sak} by using ``small'' continuous
perturbations of the curves.
\end{remark}

\section{Necessary condition}

In this section, we formulate a necessary condition for  
$C^0$ transversality for mappings of compact metric spaces into $\R^n$. We formulate our condition in terms of topological
dimension; let us recall the necessary definitions and facts (see, for example, \cite{Alexandrov} or \cite{GW}).).

\begin{definition}
Let $(X, \dist)$ be a compact metric space. The topological dimension of $X$ is the least integer $n\geq -1$
having the following property: For any $\epsilon>0$ there is an open $\epsilon$-covering of the space $X$
whose multiplicity does not exceed $n+1$.
\end{definition}

We need the following classical result of dimension theory (for example, one can find a proof in \cite{GW}, Chap. ).

\begin{proposition}
Let $(X, \dist)$ be a compact metric space and let $f\colon X \to \R^n$ be a continuous mapping of
$X$ into the $n$-dimensional Euclidean space. Then for any $\ep>0$ there exists a continuous mapping
$$
\widetilde{f}\colon X \to \R^n
$$
such that
$$
|\widetilde{f},f|_{C^0}<\epsilon,
$$
and the image $\Image\widetilde{f}$ is contained in a polyhedron $K$ whose dimension does not exceed $\dim (X)$. 
\end{proposition}

Now let $A$ and $B$ be compact metric spaces and let
$$
f\colon A \to \R^n\quad\mbox{and}\quad g\colon B \to \R^n$$
be two continuous mappings such that
$f(a)=g(b)$ for some $a\in A$ and $b\in B$.

Our necessary condition for $C^0$ transversality is as follows.

\begin{proposition}
If the mappings $f$ and $g$ are $C^0$ transverse at the pair
$(a,b)$, then for any open sets $U_A\subset A$ and $U_B\subset B$ such that $a\in U_A$ and $b\in U_B$,
$$
\dim (Cl(U_A)) + \dim (Cl(U_B)) \geq n.
$$
\end{proposition}

\textit{Proof. } Assume the converse: Let there exist neighborhoods $U_A$ and $U_B$ of the points $a$ and $b$,
respectively, such that
$$
\dim (Cl(U_A))+\dim (Cl(U_B))<n.
$$

We claim that, for any $\delta>0$, the intersection $f(U_A)\cap g(U_B)$ is not $\delta$-essential.

Let $\delta>0$. By Proposition 3, there exists a continuous mapping
$$
f_1\colon U_A \to \R^n
$$
such that
$$
|f_1,f|_{C^0}<\frac{\delta}{2}
$$
and $\Image f_1$ is contained in the union of a finite set of simplices
whose dimension does not exceed $\dim (Cl(U_A))$:
$$
f_1(U_A)\subset \cup_{i=1}^{N_A} \Delta_i^A, 
$$
where $\De_i^A=\conv \{ x_1,\dots, x_{k_i+1} \}$, $k_i\leq \dim (Cl(U_A))$.

Similarly, there exists a continuous mapping $\widetilde{g}\colon U_B \to \R^n$ such that
$$
|\widetilde{g},g|_{C^0}<\frac{\delta}{2}
\quad\mbox{and}\quad
\widetilde{g}(U_A)\subset \cup_{i=1}^{N_B} \Delta_i^B,
$$
where $\De_i^B=\conv \{ y_1,\dots, y_{l_i+1} \}$, $l_i\leq \dim (Cl(U_B))$.

Since $\dim(\De_i^A) +\dim(\De_j^B)<n$ for any $i,j$, there exists an arbitrarily small vector $v\in \R^n$ 
such that
$$
\left(  v+ \De_i^A \right) \cap \De_j^B = \varnothing
$$
for any $i,j$. If we take such a vector $v$ with $|v|<\frac{\de}{2}$, then the mapping
$$
\widetilde{f}=f_1+v
$$
satisfies the relations
$$
|\widetilde{f},f|_{C^0}\leq |\widetilde{f},f_1|_{C^0}+|f_1,f|_{C^0}<\de
$$
and
$$\Image \widetilde{f} \cap \Image \widetilde{g} = \varnothing.
$$

This proves that for any $\delta>0$, the intersection $f(U_A)\cap g(U_B)$ is not $\de$-essential.
$\Box$ 

\medskip
{\bf Aknowledgments.} The research of the authors was supported by the RFBR (project 12-01-00275)
ans by the SPbGU program ``Stability of dynamical systems with respect to perturbations and
applications to study of applied problems'' (IAS 6.38.223.2014).


\begin{thebibliography}{widestlabel}
\bibitem{Sak}
\newblock K. Sakai,
\newblock Shadowing Property and Transversality Condition,
\newblock \emph{Dynamical Systems and Chaos} (World Sci., Singapore) {\bf 1} (1995) 233–238.

\bibitem{PilSpDS}
\newblock S. Yu. Pilyugin,
\newblock ``Spaces of Dynamical Systems,''
\newblock De Gruyter, Berlin, 2012.


\bibitem{PilSak}
\newblock S. Yu. Pilyugin and K. Sakai,
\newblock $C^0$ transversality and shadowing properties,
\newblock \emph{Proceeding of the Steklov Institute of Mathematics} {\bf 256} (2007), 290-305.

\bibitem{PilSakTar}
\newblock S. Yu. Pilyugin, K. Sakai, and O. A. Tarakanov,
\newblock Transversality properties and $C^1$-open sets of diffeomorphisms with weak shadowing,
\newblock \emph{Discrete Contin. Dynam. Systems} {\bf 16} (2006), 871-882.


\bibitem{Hirsch}
\newblock Morris W. Hirsch,
\newblock ``Differential Topology,''
\newblock Springer-Verlag, New York, 1997.

\bibitem{Alexandrov}
\newblock P. S. Alexandrov and B. A. Pasynkov,
\newblock ``Introduction to the Dimension Theory,'' 
\newblock Nauka, Moscow, 1973. [in Russian]

\bibitem{GW}
\newblock W. Gurewicz and H. Wallman,
\newblock ``Dimension Theory,''
\newblock  Princeton Univ. Press, Princeton, N.J., 1941. 



\end{thebibliography}
\end{document}